\documentclass[12pt,preprint]{aastex}
\usepackage{spr-astr-addons}
\usepackage{aas_macros}
\newcommand{\abc}{{\Bigl(1+\frac{5}{2}A_2\Bigr) }}
\newcommand{\aac}{{(1-\frac{A_2}{2}) }}
\newcommand{\amc}{{(1-\mu) }}
\newcommand{\adc}{{\frac{\delta ^2}{2} }}

\newcommand{\zx}{{(x+\mu)}}
\newcommand{\zox}{{(x+\mu-1)}}
\newcommand{\zd}{{\displaystyle}}

\newcommand{\az}{{1,0}}
\newcommand{\za}{{0,1}}

\shorttitle{Nonlinear Stability ... Poynting-Robertson Drag}

\shortauthors{B.S. Kushvah et al.}

\begin{document}
\title{Nonlinear Stability in the Generalised Photogravitational Restricted Three Body Problem with Poynting-Robertson Drag}

\author{B.S.Kushvah\altaffilmark{1,2} }
\affil{Gwalior Engineering College  Airport  Road, Maharajpura,
Gwalior (M.P.)-474015,INDIA}\email{bskush@gmail.com}

\author{J.P. Sharma\altaffilmark{3} and B.Ishwar\altaffilmark{4}}
\affil{University Department of Mathematics, B.R.A. Bihar University
Muzaffarpur-842001} \email{ishwar\_bhola@hotmail.com}

\altaffiltext{1}{Ex-Junior Research Fellow Department of Science and
Technology Project} \altaffiltext{2}{Present
address for correspondence : C/O Er. Prashant Anand Kushwaha,\\
139/I, Anupam Nagar Ext-2, Opp. Jiwaji University, City Center,
Gwalior-474011} \altaffiltext{3}{Co-Principal Investigator
Department of Science and Technology Project}
\altaffiltext{4}{Principal Investigator Department of Science and
Technology Project}

%

\begin{abstract}
The Nonlinear stability of triangular equilibrium points has been
discussed in the generalised photogravitational restricted three
body problem with Poynting-Robertson drag. The problem is
generalised in the sense that smaller primary is supposed to be an
oblate spheroid. The bigger primary is considered  as radiating. We
have performed first and second order normalization of the
Hamiltonian of the problem. We have applied KAM  theorem to examine
the condition of non-linear stability. We have found three critical
mass ratios. Finally we conclude that triangular points are stable
in the nonlinear sense except three critical mass ratios at which
KAM theorem fails.
\end{abstract}



%
%
%

\section{Introduction}
The simplest form of the three-body problem is called the restricted
three-body problem(RTBP), in which a particle of infinitesimal mass
moves in the gravitational field of two massive bodies orbiting
according to the exact solution of the two-body problem. In the
circular problem, the two finite masses are fixed in a coordinate
system rotating at the orbital angular velocity, with the origin
(axis of rotation) at the centre of mass of the two bodies. Lagrange
showed that in this rotating frame there are five stationary points
at which the massless particle would remain fixed if placed there.
There are three such points lying on the line connecting the two
finite masses: one between the masses and one outside each of the
masses. The other two stationary points, called the triangular
points, are located equidistant from the two finite masses at a
distance equal to the finite mass separation they are stable in
classical case. The two masses and the triangular stationary points
are thus located at the vertices of equilateral triangles in the
plane of the circular orbit. There is a group of enthusiasts who
want to setup a colony at $L_5$ point of the Earth-Moon system. As
already noted, because $L_4$ and $L_5$ are the stable points of
equilibrium, they have been proposed for sites of large
self-contained \lq\lq Space colonies\rq\rq, an idea developed and
advocated by the late \cite{O'Neill1974}. The three body problem
have an interesting application for artificial satellites and future
space colonization. Triangular points of the Sun- Jupiter or
Sun-Earth system would be convenient sites to locate future space
colonies. Application of results to realistic actual problem is
obvious.

The classical restricted three body problem is generalized to
include the force of radiation pressure, the Poynting-Robertson(P-R)
effect and oblateness effect. The photogravitational restricted
three body problem arises from the classical problem when at least
one of the interacting bodies exerts radiation pressure, for
example, binary star systems(both primaries radiating).  The
photogravitational restricted three body problem under different
aspects was studied by \citet*{Radzievskii1950}, \citet*{Chernikov1970},
\citet*{Bhatnagar1979}, \citet*{Schuerman1980},
\citet*{KushvahBR2006},\citet*{KushvahetalNr2007}

 The Poynting-Robertson drag
named after John Henry Poynting and Howard Percy Robertson, is a
process by which solar radiation causes dust grains in a solar
system to slowly spiral inward.  \cite{Poynting1903} considered the
effect of the absorption and subsequent re-emission of sunlight by
small isolated particles in the solar system. His work was later
modified by \cite{Robertson1937} who used  precise relativistic
treatments of the first order in the ratio of the velocity of the
particle to that of light.

The location and stability of the five Lagrangian equilibrium points
in the planar, circular restricted three-body problem was
investigated by \citet*{Murray1994} when the third body is acted on
by a variety of drag forces. The approximate locations of the
displaced equilibrium points are calculated for small mass ratios
and a simple criterion for their linear stability is derived. They
showed if $a_1$ and $a_3$ denote the coefficients of the linear and
cubic terms in the characteristic equation derived from a linear
stability analysis, then an equilibrium point is asymptotically
stable provided $0<a_1< a_3$. In cases where $a_1$ is approximately
equal to $0$ or $a_1$ is approximately equal to $a_3$ the point is
unstable but there is a difference in the e-folding time scales of
the shifted $L_4$ and $L_5$ points such that the $L_4$ point, if it
exists, is less unstable than the $L_5$ point. The results are
applied to a number of general and specific drag forces. They have
shown that, contrary to intuition, certain drag forces produce
asymptotic stability of the displaced triangular equilibrium points,
$L_4$ and $L_5$. \citet*{KushvahBR2006} examined the linear
stability of triangular equilibrium points in the generalised
photogravitational restricted three body problem with
Poynting-Robertson drag and conclude that the triangular equilibrium
points are unstable due to Poynting-Robertson drag .
\citet*{KushvahetalHr2007} performed higher order normalizations in
the generalized photogravitational restricted three body problem
with Poynting-Robertson drag.

\cite{Deprit1967} investigated the nonlinear stability of triangular
points by applying Moser's modified version of Arnold's
theorem(1961). \cite{Bhatnagaretal1983} studied the effect of
perturbations on the nonlinear stability of triangular points.
\cite{BR1997} studied nonlinear stability in the generalized
restricted three body problem. His problem is generalized in the
sense that the infinitesimal body and one of the primaries have been
taken as oblate spheroid. \cite{SubbaRK1997} examined effect of
oblateness on the non-linear stability of $L_4$ in the restricted
three body problem . Hence we aim to study nonlinear stability of
triangular points in our problem.

To Examine the nonlinear stability of triangular  points we used the
KAM theorem[the work of \citet*{Kolmogorov1957} extended by \citet*{Arnold1961},
\citet*{Moser1962}]. Moser's conditions are utilised in this study by
employing the iterative scheme of Henrard for transforming the
Hamiltonian to the Birkhoff's normal form with the help of double D'
Alembert's series. We have found the second order coefficients in
the frequencies. For this we have obtained the partial differential
equations which are satisfied by the third order homogeneous
components of the fourth order part of Hamiltonian $H_4$ and second
order polynomials in the frequencies. We have found the coefficients
of sine and cosine in the homogeneous components of  order three.
They are critical terms. We have eliminated  these critical terms by
choosing properly the coefficients in the polynomials. Then we have
obtained the values of the coefficients $A,B,C$ occurring in the
fourth order part of the normalized Hamiltonian in KAM theorem. We
have applied KAM theorem to examine the conditions of nonlinear
stability. Using the first condition of the theorem, we have found
two critical mass ratios $\mu_{c1},\ \mu_{c2}$ where this condition
fails. By taking the second order coefficients, we have calculated
the determinant $D$ occurring in the second condition of the
theorem. From this, we have found the third critical mass ratio
$\mu_{c3}$ where the second condition of the theorem fails. We
conclude that triangular points are stable for all mass ratios in
the range of stability  except three critical mass ratios where KAM
theorem fails. The stability conditions are different from classical
case and others, due to radiation pressure, oblateness and P-R drag.
\section{First Order Normalization}
 \label{sec:1storder}
We used \cite{Whittaker1965} method  for the transformation of $H_2$
into the normal form

Equations of motion are as in \citet*{KushvahBR2006} and given by
\begin{eqnarray}
\ddot{x}-2n\dot{y}&=U_x ,\quad \mbox{where},\quad U_x=\frac{\partial{U_1}}{\partial{x}}-\frac{W_{1}n_1}{r^2_1}\\
\ddot{y}+2n\dot{x}&=U_y,\hspace{.85in}U_y=\frac{\partial{U_1}}{\partial{y}}-\frac{W_{1}n_2}{r^2_1}\\
U_1&=\zd{\frac{n^2(x^2+y^2)}{2}}+\frac{\amc{q_1}}{r_1}+\frac{\mu}{r_2}+\frac{\mu{A_2}}{2r^3_2}
\end{eqnarray}
\begin{eqnarray*}
r^2_1&=&\zx^2+y^2, r^2_2=\zox^2+y^2,\\\nonumber
n^2&=&1+\frac{3}{2}A_2,\\\nonumber
n_1&=&\frac{\zx[\zx\dot{x}+y\dot{y}]}{r^2_1}+\dot{x}-ny,\\
n_2&=&\frac{y[\zx\dot{x}+y\dot{y}]}{r^2_1}+\dot{y}+n\zx
\end{eqnarray*}
$W_1=\frac{(1-\mu)(1-q_1)}{c_d}$,
$\mu=\frac{m_2}{m_1+m_2}\leq\frac{1}{2}$, $m_1,m_2$ be  the masses
of the primaries, $A_2=\frac{r^2_e-r^2_p}{5r^2}$ be the oblateness
coefficient, $r_e$ and $r_p$ be the equatorial and polar radii
respectively $r$ be the distance between primaries, $c_d=
299792458$ be the dimensionless  velocity of light,
$q_1=\bigl(1-\frac{F_p}{F_g}\bigr)$ be the mass reduction factor
expressed in terms of the particle's radius $a$, density $\rho$ and
radiation pressure efficiency factor $\chi$ (in the C.G.S.system)
i.e., $q_1=1-\zd{\frac{5.6\times{10^{-5}}\chi}{a\rho}}$. Assumption
$q_1=constant$ is equivalent to neglecting fluctuation in the beam
of solar radiation, the effect of the planet's shadow, obviously
$q_1\leq1$. Triangular equilibrium points are given by
$U_x=0,U_y=0,y\neq{0}$, then we have
\begin{eqnarray}
 x_*&=&x_0\Biggl\{1\nonumber\\&&-\zd{\frac{nW_1\bigl[\amc\abc+\mu\aac\adc\bigr]}{3\mu\amc{y_0 x_0}}}\nonumber\\&&-\adc\frac{A_2}{x_0}\Biggr\} \label{eq:1x}\\
y_*&=&y_0\Biggl\{1\nonumber\\&&-\zd{\frac{nW_1\delta^2\bigl[2\mu-1-\mu(1-\frac{3A_2}{2})\adc+7\amc\frac{A_2}{2}\bigr]}{3\mu\amc{y^3_0}}}\nonumber\\&&-\zd{\frac{\delta^2\bigl(1-\adc)A_2}{y^2_0}}\Biggr\}^{1/2}\label{eq:Ly}\end{eqnarray}
where $x_0=\adc-\mu$,
$y_0=\pm\delta\bigl(1-\frac{\delta^2}{4}\bigr)^{1/2}$ and
$\delta=q^{1/3}_1$, as in  \cite{KushvahBR2006}

 The Lagrangian function of the problem can be
written as
\begin{eqnarray}
L&=&\frac{1}{2}(\dot{x}^2+\dot{y}^2)+n(x\dot{y}-\dot{x}y)+\frac{n^2}{2}(x^2+y^2)\nonumber\\&&+\frac{\amc{q_1}}{r_1}+\frac{\mu}{r_2}+\frac{\mu{A_2}}{2r^3_2}\nonumber\\
&&+W_1\Bigl\{\frac{\zx\dot{x}+y\dot{y}}{2r^2_1}-n
\arctan{\frac{y}{\zx}}\Bigr\}\nonumber\\
\end{eqnarray}
and the Hamiltonian  is $H=-L+p_x\dot{x}+p_y\dot{y}$, where
$p_x,p_y$ are the momenta coordinates given by \[
p_x=\frac{\partial{L}}{\partial{\dot{x}}}=\dot{x}-ny+\frac{W_1}{2r_1^2}\zx,
\]\[
p_y=\frac{\partial{L}}{\partial{\dot{y}}}=\dot{y}+nx+\frac{W_1}{2r_1^2}y
\]
For simplicity we suppose  $q_1=1-\epsilon$, with $|\epsilon|<<1$
then coordinates of triangular equilibrium point $L_4$  can be
written in the form
\begin{eqnarray}
x&=&\frac{\gamma}{2}-\frac{\epsilon}{3}-\frac{A_2}{2}+\frac{A_2
\epsilon}{3}\nonumber\\&&-\frac{(9+\gamma)}{6\sqrt{3}}W_1-\frac{4\gamma
\epsilon}{27\sqrt{3}}W_1 \\
y&=&\frac{\sqrt{3}}{2}\Bigl\{1-\frac{2\epsilon}{9}-\frac{A_2}{3}-\frac{2A_2
\epsilon}{9}\nonumber\\&&+\frac{(1+\gamma)}{9\sqrt{3}}W_1-\frac{4\gamma
\epsilon}{27\sqrt{3}}W_1\Bigr\}
\end{eqnarray}
where $\gamma=1-2\mu$.
 We shift the origin to $L_4$. For that, we change
$x\rightarrow {x_*}+x$ and  $y\rightarrow{y_*}+y$. Let $a=x_*+\mu,
b=y_*$ so that
\begin{eqnarray}
a&=&\frac{1}{2} \Biggl\{1-\frac{2\epsilon}{3}-A_2+\frac{2A_2
\epsilon}{3}\nonumber\\&&-\frac{(9+\gamma)}{3\sqrt{3}}W_1-\frac{8\gamma
\epsilon}{27\sqrt{3}}W_1 \bigr\}\\
b&=&\frac{\sqrt{3}}{2}\Bigl\{1-\frac{2\epsilon}{9}-\frac{A_2}{3}-\frac{2A_2
\epsilon}{9}\nonumber\\&&+\frac{(1+\gamma)}{9\sqrt{3}}W_1-\frac{4\gamma
\epsilon}{27\sqrt{3}}W_1\Bigr\}
\end{eqnarray}
Expanding $L$ in power series of $x $ and $y$, we get
\begin{eqnarray}
 L&=&L_0+L_1+L_2+L_3+\cdots \\
H&=&H_0+H_1+H_2+H_3+\cdots\nonumber\\&&
=-L+p_x{\dot{x}}+p_y{\dot{y}}
  \end{eqnarray}
  where $L_0,L_1,L_2,L_3 \ldots$ are
\begin{eqnarray}
L_0&=&\frac{3}{2}-\frac{2\epsilon}{3}-\frac{\gamma
\epsilon}{3}+\frac{ 3 \gamma A_2}{4}-\frac{3 A_2 \epsilon}{2}-\gamma
A_2 \nonumber \\
&&-\frac{\sqrt{3}W_1}{4}+\frac{2\gamma}{3\sqrt{3}}W_1\nonumber\\&&-\frac{
\epsilon W_1}{3\sqrt{3}}-\frac{23\epsilon  W_1}{54\sqrt{3}}-n
\arctan{\frac{b}{a}}
\end{eqnarray}

\begin{eqnarray}
L_1&&=\dot{x}\bigl\{-\frac{\sqrt{3}}{2}-\frac{5 A_2
}{8\sqrt{3}}+\frac{7\epsilon A_2}{12\sqrt{3}}+\frac{4
W_1}{9}-\frac{1
 \gamma W_1}{18}\bigr\}\nonumber \\&&+ \dot{y}\bigl\{\frac{1}{2}-\frac{\epsilon}{3}-\frac{A_2
}{8}+\frac{\epsilon A_2}{12\sqrt{3}}-\frac{ W_1}{6\sqrt{3}}+\frac{2
\epsilon W_1}{3\sqrt{3}}\bigr\} \nonumber \\
& &-x \bigr\{-\frac{1}{2}+\frac{\gamma}{2}+\frac{9
A_2}{8}+\frac{15\gamma A_2}{8}-\frac{35\epsilon
A_2}{12}\nonumber\\&&-\frac{29\gamma \epsilon A_2}{12}+
\frac{3\sqrt{3}W_1}{8}-\frac{2\gamma}{3\sqrt{3}}W_1-\frac{5 \epsilon
 W_1}{12\sqrt{3}}\nonumber\\&&-y
\bigr\{\frac{15\sqrt{3}A_2}{2}+\frac{9\sqrt{3}\gamma
A_2}{8}-2\sqrt{3} \epsilon A_2-2\sqrt{3}\gamma \epsilon
A_2\nonumber\\&&- \frac{W_1}{8}+\gamma W_1-\frac{43 \epsilon
}{36}W_1 \bigr\}
\end{eqnarray}
\begin{eqnarray}
L_2&=&\frac{(\dot x^2+ \dot y^2)}{2}+n(x\dot y-\dot x y)+
\frac{n^2}{2}(x^2+y^2)\nonumber\\&&-Ex^2-Fy^2-G xy
\end{eqnarray}
\begin{equation}
L_3=-\frac{1}{3!}\left\{x^3T_1+3x^2yT_2+3xy^2T_3+y^3T_4+6T_5\right\}\label{eq:L3}
  \end{equation}
  \begin{eqnarray}
L_4&=&-\frac{1}{4!}\Bigr\{N_1x^4+4N_2x^3y+6N_3x^2y^2\nonumber\\&&+4N_4xy^3+24N_6\Bigr\}\label{eq:L4}
  \end{eqnarray}
  where
\begin{eqnarray}
E&=&\frac{1}{16}\Bigl[ 2-6\epsilon- 3A_2-
\frac{31A_2\epsilon}{2}-\frac{69W_1}{6\sqrt{3}} \nonumber\\&&+\gamma
\bigl\{2\epsilon+12A_2+
\frac{A_2\epsilon}{3}+\frac{199W_1}{6\sqrt{3}}\bigr\}\Bigr]
  \end{eqnarray}
  \begin{eqnarray}
F&=&\frac{-1}{16}\Bigl[ 10-2\epsilon+21A_2-
\frac{717A_2\epsilon}{18}-\frac{67W_1}{6\sqrt{3}} \nonumber
\\&+&\gamma \bigl\{6\epsilon-
\frac{293A_2\epsilon}{18}+\frac{187W_1}{6\sqrt{3}}\bigr\}\Bigr]
  \end{eqnarray}
   \begin{eqnarray}
G&=&\frac{\sqrt{3}}{8}\Bigl[2\epsilon+6A_2-
\frac{37A_2\epsilon}{2}-\frac{13W_1}{2\sqrt{3}} \nonumber
\\&&-\gamma \bigl\{6\epsilon-\frac{\epsilon}{3}+13A_2-
\frac{33A_2\epsilon}{2}+\frac{(11W_1}{2\sqrt{3}}\bigr\}\Bigr]\nonumber
\\
  \end{eqnarray}
$T_i,N_j,(i=1,\dots,5,\ j=1,\dots,6 )$ are as in
\citet*{KushvahetalHr2007}.

The second order part $H_2$ of the corresponding Hamiltonian takes
the form
\begin{equation}
H_2=\frac{p_x^2+p_y^2}{2}+n(yp_x-xp_y)+Ex^2+Fy^2+Gxy
\end{equation}
To investigate the stability of the motion, as in
\cite{Whittaker1965}, we consider the following set of linear
equations in the variables $x, y$:
 \begin{equation*}
 \begin{array}{l c l}
 -\lambda p_x& = & \frac{\partial{H_2}}{\partial x},\\&&\\
 -\lambda p_y& = & \frac{\partial{H_2}}{\partial y},\\
 \end{array} \quad
 \begin{array}{l c l }
 \lambda x& = & \frac{\partial{H_2}}{\partial p_x}\\&&\\
 \lambda y& = & \frac{\partial{H_2}}{\partial p_y}\\
 \end{array}
 \end{equation*}
 i.e.
 \begin{equation} AX=\mathbf{0} \label{eq:ax}\end{equation}
where
  \begin{equation}
  X=\left[\begin{array}{c}
  x\\
  y\\
  p_x\\
  p_y \end{array}\right], \quad
 A=\left[\begin{array}{c c c c}
  2E & G&\lambda& -n\\
G&2F&n&\lambda\\
  -\lambda& n& 1& 0\\
  -n & -\lambda& 0& 1\end{array}\right]
  \end{equation}

 Clearly $|A|=0$, implies that the characteristic equation
 corresponding to Hamiltonian $H_2$ is given by
 \begin{equation}
 \lambda^4+2(E+F+n^2)\lambda^2+4EF -G^2+n^4-2n^2(E+F)=0 \label{eq:ch}
 \end{equation}
 This is characteristic equation whose discriminant is
  \begin{equation}
 D=4(E+F+n^2)^2-4\bigl\{4EF-G^2+n^4-2n^2(E+F)\bigr\}
 \end{equation}
 Stability is assured  only when $D>0$.
 i.e
  \begin{eqnarray}
  \mu<\mu_{c_0}-0.221896\epsilon +2.103887A_2
  +\nonumber\\
    0.493433\epsilon A_2 +0.704139 W_1 +
    0.401154\epsilon W_1\label{eq:muc0}
 \end{eqnarray}
 where $ \mu_{c_0}=0.038521$,(Routh's
critical mass ratio)
 When $D>0$ the roots $\pm i \omega_1$ and $\pm i \omega_2$ ($\omega_1,$ $\omega_2$ being the long/short -periodic frequencies) are related to each other  as

 \begin{eqnarray}
   \omega_1^2+\omega_2^2&=& 1-\frac{\gamma \epsilon}{2}+\frac{3\gamma A_2}{2}+\frac{83\epsilon A_2}{12}-\frac{W_1}{24\sqrt{3}}\nonumber\\
  \label{eq:w1+w2}
 \end{eqnarray}
 \begin{eqnarray}
  \omega_1^2\omega_2^2&=&\frac{27}{16} -\frac{27\gamma^2}{16}+\frac{9\epsilon}{8}+\frac{9\gamma\epsilon}{8} +\frac{117\gamma A_2}{16}\nonumber\\
  &&-\frac{241\epsilon A_2}{32}+\frac{35W_1}{16\sqrt{3}}-\frac{55 \sqrt{3}\gamma W_1}{16}\label{eq:w1w2}\\
&&(0<\omega_2<\frac{1}{\sqrt{2}}<\omega_1<1)\nonumber\end{eqnarray}
From (~\ref{eq:w1+w2}) and  (~\ref{eq:w1w2}) it may be noted that
$\omega_j$ $ (j=1,2)$ satisfy
\begin{eqnarray}
  &&\gamma^2= 1+\frac{4\epsilon}{9}-\frac{107\epsilon A_2}{27}+\frac{2\gamma \epsilon }{3}-\frac{25W_1}{27\sqrt{3}}\nonumber\\&&
  +\biggl(-\frac{16}{27}+\frac{32\epsilon}{243}+\frac{208 A_2}{81}\nonumber\\&&-\frac{8\gamma A_2}{27}-\frac{4868\epsilon A_2}{729}+\frac{296W_1}{243\sqrt{3}}\biggr)\omega_j^2\nonumber\\
  && +\biggl(\frac{16}{27}-\frac{32\epsilon}{243}-\frac{208 A_2}{81}\nonumber\\&&-\frac{1880\epsilon A_2}{729}-\frac{2720W_1}{2187\sqrt{3}}\biggr)\omega_j^4
 \end{eqnarray}
 Alternatively, it can also be seen that if $u=\omega_1\omega_2$,
 then (~\ref{eq:w1w2}) gives
\begin{eqnarray}
  \gamma^2&=& 1+\frac{4\epsilon}{9}-\frac{107\epsilon A_2}{27}-\frac{25W_1}{27\sqrt{3}}\nonumber\\&&+\gamma\biggl(\frac{2\epsilon }{3}+\frac{1579\epsilon A_2}{324}-\frac{55\gamma W_1}{9\sqrt{3}}\biggr)\nonumber\\&&
  +\biggl(-\frac{16}{27}+\frac{32\epsilon}{243}+\frac{208 A_2}{81}\nonumber\\&&-\frac{1880\epsilon A_2}{729}+\frac{320W_1}{243\sqrt{3}}\biggr)u^2 \end{eqnarray}
Following the method for reducing $H_2$ to the normal form, as in
\cite{Whittaker1965},use the transformation
\begin{equation}X=JT \label{eq:XJT}\end{equation}
\begin{equation*}
\mbox{where}\, X=\left[\begin{array}{c}
x\\y\\p_x\\p_y\end{array}\right],J=[J_{ij}]_{1\leq i, j \leq 4},\
T=\left[\begin{array}{c} Q_1\\Q_2\\P_1\\P_2\end{array}\right]
\end{equation*}
where $J_{ij}$ are as in \citet*{KushvahetalHr2007},\, $P_i= (2
I_i\omega_i)^{1/2}\cos{\phi_i},$ $Q_i= (\frac{2
I_i}{\omega_i})^{1/2}\sin{\phi_i},$ $(i=1,2)$

The transformation changes the second order part of the Hamiltonian
into the normal form \begin{equation}
H_2=\omega_1I_1-\omega_2I_2\end{equation} The general solution of
the corresponding equations of motion are
\begin{equation} I_i=\mbox{const.}, \quad \phi_i=\pm \omega_i+\mbox{const.},\  (i=1,2) \label{eq:intmt} \end{equation}
If the oscillations about $L_4$ are exactly linear, the
Eq.(~\ref{eq:intmt}) represent the integrals of motion and the
corresponding orbits will be given by
\begin{eqnarray}x&=&J_{13}\sqrt{2\omega_1I_1}\cos{\phi_1}+J_{14}\sqrt{2\omega_2I_2}\cos{\phi_2}\label{eq:xb110}\end{eqnarray}
\begin{eqnarray}
y&=&J_{21}\sqrt{\frac{2I_1}{\omega_1}}\sin{\phi_1}+J_{22}\sqrt{\frac{2I_2}{\omega_2}}\sin{\phi_2}\nonumber\\&&+J_{23}\sqrt{2I_1}{\omega_1}\cos{\phi_1}+J_{24}\sqrt{2I_2}{\omega_2}\sin{\phi_2}\label{eq:yb101}
\end{eqnarray}
\section{Second Order Normalization}
\label{sec:2ndorder}
  In order to perform Birkhoff's normalization, we use Henrard's
  method (\citet{Deprit1967}) for which the
  coordinates $(x,y)$ of infinitesimal body, to be expanded in
  double D'Alembert series  $x=\sum_{n\geq1}B_n^\az,\quad y=\sum_{n\geq 1}B_n^\za$
  where the homogeneous components $B_n^\az$ and $B_n^\za$ of
  degree $n$ are of the form
  \begin{eqnarray}
&&\sum_{0\leq{m}\leq{n}}
I_1^{\frac{n-m}{2}}I_2^{\frac{m}{2}}\sum_{(p,q)}\bigl[C_{n-m,m,p,q}
\cos{(p\phi_1+q\phi_2)}\nonumber\\&&+S_{n-m,m,p,q}
\sin{(p\phi_1+q\phi_2)}\bigr]
  \end{eqnarray}
  The conditions in double summation are (i) $p$ runs over those
  integers in the interval $0\leq p\leq n-m$ that have the same
  parity as $n-m$ (ii) $q$ runs over those integers in the interval $-m\leq q\leq
  m$ that have the same parity as $m$. Here $I_1$, $I_2$ are the
  action momenta coordinates which are to be taken as constants of
  integer, $\phi_1$, $\phi_2$ are angle coordinates to be
  determined as linear functions of time in such a way that $\dot\phi_1=\omega_1+\sum_{n\geq 1}f_{2n}(I_1,I_2),\ \dot\phi_2=-\omega_2+\sum_{n\geq 1}g_{2n}(I_1,I_2)$  where  $\omega_1,\omega_2$ are the basic  frequencies, $f_{2n}$ and  $g_{2n}$ are of the form
  \begin{eqnarray}
  f_{2n}&=&\sum_{0\leq m\leq n}{f'}_{2(n-m),2m}I_1^{n-m}I_2^m\\
  g_{2n}&=&\sum_{0\leq m\leq n}{g'}_{2(n-m),2m}I_1^{n-m}I_2^m
  \end{eqnarray}
The first order components $B_1^\az$ and $B_1^\za$ are the values of
 $x$ and  $y$ given by (~\ref{eq:xb110}) (~\ref{eq:yb101}). In order to find out  the second order components $B_2^\az,B_2^\za$
we consider  Lagrange's  equations of motion
\begin{equation}
\frac{d}{dt}(\frac{\partial L}{\partial \dot x })-\frac{\partial
L}{\partial x }=0, \quad \frac{d}{dt}(\frac{\partial L}{\partial
\dot y })-\frac{\partial L}{\partial y }=0 \end{equation}
\begin{equation}
{i.e.}\,\left.\begin{array}{l c l}
\ddot x-2n\dot y+(2E-n^2)x+Gy&=&\frac{\partial L_3}{\partial x }+\frac{\partial L_4}{\partial x }\\
&&\\\ddot x+2n\dot x+(2F-n^2)y+Gx&=&\frac{\partial L_3}{\partial y
}+\frac{\partial L_4}{\partial y }\end{array}
\right\}\label{eq:lgeq}\end{equation} Since $x$ and $y$ are double
D'Alembert series, $x^jx^k(j\geq0,k\geq0,j+k\geq0)$ and  the time
derivatives $\dot x ,\dot y ,\ddot x, \ddot y $ are also double
D'Alembert series. We can write
\[\dot x=\sum_{n\geq 1} \dot x_n,\,\dot y=\sum_{n\geq 1} \dot
y_n,\,\ddot x=\sum_{n\geq 1} \ddot x_n,\,\ddot y=\sum_{n\geq 1}
\ddot y_n \] where $\dot x ,\dot y ,\ddot x, \ddot y $ are
homogeneous components of degree $n$ in $I_1^{1/2},I_2^{1/2}$ i.e.
\begin{eqnarray} \dot x &=&
\frac{d}{dt}\sum_{n\geq 1}B_n^\az\nonumber\\&=&\sum_{n\geq
1}\Biggl[\frac{\partial
B_n^\az}{\partial{\phi_1}}(\omega_1+f_2+f_4+\cdots)\nonumber\\&+&\frac{\partial
B_n^\az}{\partial{\phi_2}}(-\omega_2+g_2+g_4+\cdots)\Biggr]\end{eqnarray}
We write three components $\dot x_1 ,\dot x_2 ,\dot x_3$ of $\dot x$
\begin{eqnarray}
\dot x_1&=&\omega_1\frac{\partial
B_1^\az}{\partial{\phi_1}}-\omega_2\frac{\partial
B_1^\az}{\partial{\phi_2}}=DB_1^\az\\
 \dot
x_2&=&\omega_1\frac{\partial
B_2^\az}{\partial{\phi_1}}-\omega_2\frac{\partial
B_2^\az}{\partial{\phi_2}}=DB_2^\az\\
\dot x_3&=&\omega_1\frac{\partial
B_3^\az}{\partial{\phi_1}}-\omega_2\frac{\partial
B_3^\az}{\partial{\phi_2}}+f_2\frac{\partial
B_1^\az}{\partial{\phi_1}}-g_2\frac{\partial
B_1^\az}{\partial{\phi_2}}\nonumber\\
&=&DB_2^\az+f_2\frac{\partial
B_1^\az}{\partial{\phi_1}}-g_2\frac{\partial
B_1^\az}{\partial{\phi_2}}
\end{eqnarray}
where \begin{equation}D\equiv \omega_1\frac{\partial\
}{\partial{\phi_1}}-\omega_2\frac{\partial\
}{\partial{\phi_2}}\end{equation}
 Similarly three components $\ddot
x_1 ,\ddot x_2 ,\ddot x_3$ of $\ddot x$ are
\begin{eqnarray*}
\ddot x_1 &=&D^2B_1^\az,\, \ddot x_2=D^2B_2^\az,\\&& \ddot
x_3=D^2B_3^\az+2\omega_1f_2\frac{\partial^2B_1^\az}{\partial\phi_1^2}-2\omega_2g_2\frac{\partial^2B_1^\az}{\partial\phi_2^2}
\end{eqnarray*}
In similar manner we can write the components of $\dot y, \ddot y$.
Putting the values of  $x, y, \dot x ,\dot y ,\ddot x $ and $\ddot
y$ in terms of double D'Alembert series in Eq.(~\ref{eq:lgeq}) we
get
\begin{eqnarray}
&&\left(D^2+2E-1-\frac{3}{2}A_2\right)B_2^\az\nonumber\\&&-\left\{2\left(
1+\frac{3}{4}A_2\right)D-G\right\}B_2^\za =X_2 \label{eq:x2}
\end{eqnarray}
\begin{eqnarray}&&
\left\{2\left(
1+\frac{3}{4}A_2\right)D+G\right\}B_2^\az\nonumber\\&&+\left(D^2+2F-1-\frac{3}{2}A_2\right)B_2^\za=Y_2\label{eq:y2}
\end{eqnarray}

 where \[X_2=\left[\frac{\partial
L_3}{\partial x}\right]_{x=B_1^\az,y=B_1^\za},\,\,
Y_2=\left[\frac{\partial L_3}{\partial
y}\right]_{x=B_1^\az,y=B_1^\za}\] These are two simultaneous partial
differential equations in $B_2^\az$ and $B_2^\za$. We solve these
equations to find the values of $B_2^\az$ and $B_2^\za$, from
(~\ref{eq:x2}) and (~\ref{eq:y2})
\begin{equation}
\triangle_1 \triangle_2B_2^\az=\Phi_2,\, \triangle_1
\triangle_2B_2^\za=-\Psi_2 \label{eq:phi_si}\end{equation}
\[\mbox{where} \quad \triangle_1=D^2+\omega_1^2,
\triangle_2=D^2+\omega_2^2
\]
\begin{equation}
\Phi_2=(D^2+2F-n^2)X_2+(2nD-G)Y_2 \label{eq:phi2}
\end{equation}
\begin{equation}
\Psi_2=(2nD+G)X_2-(D^2+2E-n^2)Y_2 \label{eq:psi2}
\end{equation}
The Eq.(~\ref{eq:phi_si}) can  be solved for $B_2^\az$ and $B_2^\za$
by putting the formula
\begin{equation*}\frac{1}{\triangle_1\triangle_2}\left\{\begin{array}{c}\cos(p\phi_1+q\phi_2)\\ \mbox{or} \\\sin(p\phi_1+q\phi_2)\end{array}=\frac{1}{\triangle_{p,q}}\left\{\begin{array}{c}\cos(p\phi_1+q\phi_2)\\\mbox{or} \\\sin(p\phi_1+q\phi_2)\end{array}\right.\right.\]
where \[\triangle_{p,q}=\left[
\omega_1^2-(\omega_1p-\omega_2q)^2\right]\left[
\omega_2^2-(\omega_1p-\omega_2q)^2\right]
\end{equation*}
provided $\triangle_{p,q}\neq0$. Since $\triangle_{1,0}=0,
\triangle_{0,1}=0$ the terms
$\cos\phi_1,\sin\phi_1,\cos\phi_2,\sin\phi_2$ are the critical
terms. $\Phi_2$ and $\Psi_2$ are free from such terms. By
condition(1) of Moser's theorem $k_1\omega_1+k_2\omega_2\neq 0$  for
all pairs $(k_1,k_2)$ of integers such that $|k_1|+|k_2|\leq4$,
therefore each of $\omega_1, \omega_2,
\omega_1\pm2\omega_2,\omega_2\pm2\omega_1$ is different from zero
and consequently none of the divisors $\triangle_{0,0},
\triangle_{0,2}, \triangle_{2,0}, \triangle_{1,1}, \triangle_{1,-1}$
is zero. The second order components $B_2^\az, B_2^\za$ are as
follows:
\begin{eqnarray}
B_2^\az&=&r_1I_1+r_2I_2+r_3I_1\cos2\phi_1\nonumber\\&&+r_4I_2\cos2\phi_2
+r_5I_1^{1/2}I_2^{1/2}\cos(\phi_1-\phi_2)\nonumber\\&&+r_6I_1^{1/2}I_2^{1/2}\cos(\phi_1+\phi_2)+r_7I_1\sin2\phi_1\nonumber\\&&+r_8I_2\sin2\phi_2
+r_9I_1^{1/2}I_2^{1/2}\sin(\phi_1-\phi_2)\nonumber\\&&+r_{10}I_1^{1/2}I_2^{1/2}\sin(\phi_1+\phi_2)\label{eq:b2az}
\end{eqnarray} \begin{eqnarray}
B_2^\za&=&-\Bigl\{s_1I_1+s_2I_2+s_3I_1\cos2\phi_1\nonumber\\&&+s_4I_2\cos2\phi_2+s_5I_1^{1/2}I_2^{1/2}\cos(\phi_1-\phi_2)\nonumber\\
&&
+s_6I_1^{1/2}I_2^{1/2}\cos(\phi_1+\phi_2)+s_7I_1\sin2\phi_1\nonumber\\&&+s_8I_2\sin2\phi_2+s_9I_1^{1/2}I_2^{1/2}\sin(\phi_1-\phi_2)\nonumber\\&&+s_{10}I_1^{1/2}I_2^{1/2}\sin(\phi_1+\phi_2)\Bigr\}\label{eq:b2za}
\end{eqnarray}
where $r_i,s_i,(i=1,\dots,10)$ are as in \citet*{KushvahetalHr2007}.
Using transformation $x=B_1^{\az}+B_2^{\az}$ and
$y=B_1^{\za}+B_2^{\za}$ the third order part $H_3$ of the
Hamiltonian in $I_1^{1/2},I_2^{1/2}$ is  of the form
\begin{equation}\label{eq:H3} H_3=A_{3,0}I_1^{3/2}+A_{2,1}I_1I_2^{1/2}+A_{1,2}I_1^{1/2}I_2+A_{0,3}I_2^{3/2}
\end{equation}
We can verify that in Eq.(~\ref{eq:H3}) $A_{3,0}$ vanishes
independently as in \cite{Deprit1967}. Similarly the other
coefficients $A_{2,1},$ $A_{1,2},$ $A_{0,3}$ are also found to be
zero independently.
\section{Second Order Coefficients in the Frequencies}
\label{sec:2coef} In order to find out the second order coefficients
$f_{2,0}, f_{0,2}, g_{2,0}, g_{0,2}$  in the polynomials $f_2$ and
$g_2$ we have done as in \cite{Deprit1967}. Proceeding as
(~\ref{eq:phi_si}), we  find
\begin{equation}
\triangle_1 \triangle_2B_3^\az=\Phi_3-2f_2P-2g_2Q\label{eq:phi3}
\end{equation}
\begin{equation}
 \triangle_1\triangle_2B_3^\za=\Psi_3-2f_2U-2g_2V \label{eq:psi3}
 \end{equation}
where
\begin{equation}
\Phi_3=\Bigl[D^2+2F-n^2\Bigr]X_3+\Bigl[(2nD-G)\Bigr]Y_3\label{eq:Phi3}\\
\end{equation}
\begin{equation}
\Psi_3=-\Bigl[2(nD+G)\Bigr]X_3+\Bigl[D^2+2nE-n^2\Bigr]Y_3\label{eq:Psi3}\end{equation}
\begin{eqnarray}
P&=&\left[ D^2+2F-n^2
\right]\left[\omega_1\frac{\partial^2B_1^\az}{\partial\phi_1^2}-n\frac{\partial
B_1^\za}{\partial
\phi_1}\right]\nonumber\\&&+(2nD-G)\left[\omega_1\frac{\partial^2B_1^\za}{\partial\phi_1^2}+n\frac{\partial
B_1^\az}{\partial \phi_1}\right]\label{eq:P}
\end{eqnarray}
\begin{eqnarray}
Q&=&-\left[ D^2+2F-n^2
\right]\left[\omega_2\frac{\partial^2B_1^\az}{\partial\phi_2^2}-n\frac{\partial
B_1^\za}{\partial
\phi_1}\right]\nonumber\\&&-(2nD-G)\left[\omega_2\frac{\partial^2B_1^\za}{\partial\phi_2^2}+n\frac{\partial
B_1^\az}{\partial \phi_2}\right]\label{eq:Q}
\end{eqnarray}
\begin{eqnarray}
U&=&-(2nD+G)\left[\omega_1\frac{\partial^2B_1^\az}{\partial\phi_1^2}-n\frac{\partial
B_1^\za}{\partial \phi_1}\right]\nonumber\\&+&\left[ D^2+2E-n^2
\right]\left[\omega_1\frac{\partial^2B_1^\za}{\partial\phi_1^2}+n\frac{\partial
B_1^\az}{\partial \phi_1}\right]\label{eq:U}
\end{eqnarray}
\begin{eqnarray}
V&=&(2nD+G)\left[\omega_2\frac{\partial^2B_1^\az}{\partial\phi_2^2}-n\frac{\partial
B_1^\za}{\partial \phi_2}\right]\nonumber\\&-&\left[ D^2+2E-n^2
\right]\left[\omega_2\frac{\partial^2B_1^\za}{\partial\phi_2^2}-n\frac{\partial
B_1^\az}{\partial \phi_2}\right]\label{eq:V}
\end{eqnarray}
\begin{equation}
X_3=\frac{\partial}{\partial x}(L_3+L_4), \quad
Y_3=\frac{\partial}{\partial y}(L_3+L_4)
\end{equation}
i.e.
\begin{eqnarray}
X_3&=&\frac{T_1}{2}x^2+T_2 xy+\frac{T_3}{2}y^2+\frac{N_1}{6}x^3+\frac{N_2}{2}x^2y\nonumber\\&&+\frac{N_3}{2}xy^2+\frac{N_4}{6}y^3+\frac{\partial{T_5}}{\partial{x}}+\frac{\partial{N_6}}{\partial{x}}\\
Y_3&=&\frac{T_2}{2}x^2+T_3
xy+\frac{T_4}{2}y^2+\frac{N_2}{6}x^3+\frac{N_3}{2}x^2y\nonumber\\&&+\frac{N_4}{2}xy^2+\frac{N_5}{6}y^3+\frac{\partial{T_5}}{\partial{y}}+\frac{\partial{N_6}}{\partial{y}}
\end{eqnarray}
(\ref{eq:Phi3}) and (\ref{eq:Psi3}) are the partial differential
equations which are satisfied by the third order components
$B_3^\az, B_3^\za$ and the second order polynomials $f_2,g_2$ in the
frequencies. We do not require to find out the components $B_3^\az$
and  $B_3^\za$. We  find the coefficients of $\cos\phi_1, \sin
\phi_1, \cos\phi_2$ and $\sin\phi_2$ in the right hand sides of
(~\ref{eq:Phi3}),(\ref{eq:Psi3}). They are the critical terms ,
since $\triangle_{1,0}=\triangle_{0,1}=0$. We eliminate these terms
by choosing properly the coefficients in the polynomials
\begin{equation} f_2=f_{2,0}I_1+f_{0,2}I_2, \quad
g_2=g_{2,0}I_1+g_{0,2}I_2\label{eq:f2}
\end{equation}
Further, we find that
\begin{eqnarray}
f_{2,0}=\frac{(\mbox{coefficient of } \cos\phi_1 \mbox { in } \Phi_3)}{2 (\mbox{ coefficient of } \cos\phi_1 \mbox{ in }P)}=A \\
f_{0,2}=g_{2,0}=\frac{(\mbox{coefficient of } \cos\phi_2 \mbox { in }
\Phi_3)}{2( \mbox{ coefficient of } \cos\phi_2 \mbox{ in } Q)}=B
\\g_{0,2}=\frac{(\mbox{coefficient of } \cos\phi_2 \mbox { in
}\Psi_3)}{2 (\mbox{ coefficient of } \cos\phi_2 \mbox{ in } Q)}=C
\\\nonumber
\end{eqnarray}
where
\begin{eqnarray}
A&=&A_{1,1}+(A_{1,2}+A_{1,3}\gamma)\epsilon+(A_{1,4}\nonumber\\&&+A_{1,5}\gamma)A_2+(A_{1,6}+A_{1,7}\gamma)W_1
\end{eqnarray}
\begin{eqnarray}
B&=&B_{1,1}+(B_{1,2}+B_{1,3}\gamma)\epsilon+(B_{1,4}\nonumber\\&&+B_{1,5}\gamma)A_2+(B_{1,6}+B_{1,7}\gamma)W_1
\end{eqnarray}
\begin{eqnarray}
C&=&C_{1,1}+(C_{1,2}+C_{1,3}\gamma)\epsilon+(C_{1,4}\nonumber\\&&+C_{1,5}\gamma)A_2+(C_{1,6}+C_{1,7}\gamma)W_1
\end{eqnarray}
where $A_{1,i},B_{1,i}$ and $C_{1,i},(i=1\dots,7)$ are as in
Appendix I
\section{Stability}
\label{sec:stab} The condition(i) of KAM theorem fails when
$\omega_1=2\omega_2$ and $\omega_1=3\omega_2$
\subsection{Case(i)}\label{subsect:i}
\begin{equation} \mbox{When}\quad \omega_1=2\omega_2\label{eq:KAMi}\end{equation}
Then from (~\ref{eq:KAMi}) and (~\ref{eq:w1w2})  we have
\begin{eqnarray}
&&\mu^2\left(-\frac{27}{4}-\frac{3\epsilon}{2}-\frac{117A_2}{4}-\frac{221W_1}{15\sqrt{3}}\right)\nonumber\\&&+\mu\left(\frac{27}{4}-\frac{107\epsilon}{100}+\frac{3021A_2}{100}+\frac{4291W_1}{120\sqrt{3}}\right)\nonumber\\&&-\frac{4}{25}+\frac{407\epsilon}{200}-\frac{12A_2}{25}-\frac{23991W_1}{200\sqrt{3}}=0
\end{eqnarray}
Solving for $\mu$ we have
\begin{eqnarray}
&&\mu_{c1}=
 0.024294 -
    0.312692\epsilon\nonumber\\&&- 0.036851A_2 + 1.001052
    W_1\label{eq:muc1}
\end{eqnarray}

\subsection{Case(ii)}\label{subsect:ii}
\begin{equation} \mbox{When} \quad\omega_1=3\omega_2\label{KAMii}\end{equation}
Proceeding as (~\ref{subsect:i}), we have
\begin{eqnarray}
&&\mu^2\left(-\frac{27}{4}-\frac{3\epsilon}{2}-\frac{117A_2}{4}-\frac{99
\sqrt{3}W_1}{20}\right)\nonumber\\&&+\mu\left(\frac{27}{4}-\frac{93\epsilon}{100}+\frac{2979A_2}{100}+\frac{119
\sqrt{3}W_1}{10}\right)\nonumber\\&&-\frac{9}{100}+\frac{393\epsilon}{200}-\frac{27A_2}{100}-\frac{4777W_1}{400\sqrt{3}}=0
\end{eqnarray}
Solving for $\mu$, we have
\begin{eqnarray}
&&\mu_{c2}=0.013516 -
    0.29724\epsilon\nonumber\\&&- 0.019383 A_2 + 1.007682 W_1\label{eq:muc2}
\end{eqnarray}
Normalized Hamiltonian up to fourth order is
\begin{eqnarray}\label{eq:H_Norm4}
H=\omega_1I_1-\omega_2I_2+\frac{1}{2}(AI_1^2+2BI_1I_2+CI_2^2)+\dots
\end{eqnarray}
Calculating the determinant $ D$ occurring in condition (ii) of KAM
theorem, we have \[D=-(A\omega_2^2+2B\omega_1\omega_2+C\omega_1^2)\]
Putting the values of $A, B$ and $C$ and if $u=\omega_1\omega_2$, we
have
\begin{eqnarray}
&&D=\frac{644u^4-541u^2+36}{8(4u^2-1)(25u^2-4)}+(D_2+D_3\gamma)\epsilon\nonumber\\&&+(D_4+D_5\gamma)A_2+(D_6+D_7\gamma)W_1\label{eq:D}
\end{eqnarray}
The second condition of KAM theorem is satisfied if , in the
interval $0<\mu<\mu_{c0}$, [where $\mu_{c0}$ as in ~(\ref{eq:muc0})]
the mass parameter does not take the value $\mu_{c3}$, which makes
$D=0$. To find $\mu_{c3}$, we note that when $\epsilon=A_2=W_1=0$,
then from (~\ref{eq:D}), $D$ becomes zero if and only if
\[644u^4-541u^2+36=0\] This implies that
\begin{equation}
u^2=\frac{541-\sqrt {199945}}{1288}=0.072863
=u_0\mbox{(say)}\end{equation} Writing
$u^2=\frac{27(1-\gamma^2)}{16}$,\, $\gamma=1-2\mu $ and then solving
above, we get
\begin{equation}
\gamma=\gamma_0=0.978173\dots,\, \mu=\mu_0=0.010914\dots
\end{equation}
When $\epsilon,A_2, W_1$ are not zero, we assume that $D$ is zero if
\begin{equation}
\mu=\mu_0+\alpha_1 \epsilon+\alpha_2A_2+\alpha_3W_1
\end{equation}
\begin{equation}
\gamma=\gamma_0-2(\alpha_1 \epsilon+\alpha_2A_2+\alpha_3W_1)
\end{equation}
\[\mbox{where}\quad \gamma_0=1-2\mu_0 \]
\begin{eqnarray}
&&u^2=u_0+(u_1+\alpha_1
u_2)\epsilon\nonumber\\&&+(u_3+\alpha_2u_4)A_2+(u_5+\alpha_3u_6)W_1\label{eq:u}
\end{eqnarray}
with \begin{eqnarray*}
&&u_1=\frac{27}{16}\gamma_0^2+\frac{9}{8}\gamma_0+\frac{9}{8},u_3=\frac{117(1-\gamma^2_0)}{16},\nonumber\\&&u_2=u_4=\frac{27\gamma_0}{4},u_6=\frac{27\gamma_0}{4\sqrt{3}},\nonumber\\&&u_5=\frac{27\gamma_0^2+165\gamma_0+35}{16\sqrt{3}}
\end{eqnarray*}
and $\alpha_i,(i=1,2,3)$ are to be determined. From (~\ref{eq:D}),
$D$ is zero when
\begin{eqnarray}
&D&=\frac{644u^4-541u^2+36}{8(4u^2-1)(25u^2-4)}+(D_2+D_3\gamma)\epsilon\nonumber\\&&+(D_4+D_5\gamma)A_2+(D_6+D_7\gamma)W_1=0\label{eq:Dzero}\end{eqnarray}
Making use of (~\ref{eq:u}) in (~\ref{eq:Dzero}) and equating to
zero the coefficients of $\epsilon,A_2$ and $W_1$, we get
\begin{eqnarray}
\alpha_1&=&-\frac{1}{u_2(1288u_0-541)}\Bigl\{(1288u_0-541)u_1\nonumber\\&&+8(D_2^0+D_3^0\gamma_0)(4u_0-1)(25u_0-4)\Bigr\}
\end{eqnarray}
\begin{eqnarray}\alpha_2&=&-\frac{1}{u_4(1288u_0-541)}\Bigl\{(1288u_0-541)u_3\nonumber\\&&+8(D_4^0+D_5^0\gamma_0)(4u_0-1)(25u_0-4)\Bigr\}
\end{eqnarray}
\begin{eqnarray}\alpha_3&=&-\frac{1}{u_6(1288u_0-541)}\Bigl\{(1288u_0-541)u_5\nonumber\\&&+8(D_6^0+D_7^0\gamma_0)(4u_0-1)(25u_0-4)\Bigr\}\end{eqnarray}
where $D_n^0$ $(n=2,3,4,\dots ,7)$ are $D_n$ given as in Appendix
II, as evaluated for the unperturbed problem. Numerical computation
yields,
\[\alpha_1=-0.120489\dots,\ \alpha_2=-0.373118\dots,\]\[
\alpha_3=2.904291\dots\] Then we  have
\begin{eqnarray}
&&\mu_{c3}=\mu_0+\alpha_1\epsilon+\alpha_2A_2+\alpha_3W_1
\nonumber\\&&=0.010914-0.120489
\epsilon\nonumber\\&&-0.373118A_2+2.904291
W_1\label{eq:muc3}\end{eqnarray} Hence in the interval
$0<\mu<\mu_{c0}$, both the conditions of KAM theorem are satisfied
and therefore the triangular point is stable except for three mass
ratios $\mu_{ci}(i=1,2,3)$.
\section{Analytical Study}\label{sec:ASob}
\subsection{Observation I}\label{subsect:obs1} Consider
$A_2=0$,$q_1=1,(W_1=0)$ then problem reduced to the classical
restricted three body problem. From equation (~\ref{eq:1x})
(~\ref{eq:Ly}) we get
 \[\quad x=\frac{1}{2}-\mu, \quad
y=\pm\frac{\sqrt{3}}{2} \] from (~\ref{eq:muc0}) stability is
assured  when $\mu<\mu_{c0}$ where $\mu_{c0}=0.038521$. The
 relation between $\omega_1,
\omega_2$ in (~\ref{eq:w1+w2})(~\ref{eq:w1w2}) are given by
\begin{eqnarray}&& \omega_1^2+\omega_2^2=1,\, \omega_1^2\omega_2^2=\frac{27}{16}{1-\gamma^2}
\\&&(0<\omega_2<\frac{1}{\sqrt{2}}<\omega_1<1)\nonumber\end{eqnarray}

From (~\ref{eq:muc1}), (~\ref{eq:muc2}) (~\ref{eq:muc3}) we have
found that the  triangular points are stable in the range of linear
stability except the three mass ratios
\begin{eqnarray}\mu_{c1}&=& 0.024294\label{eq:muc1class} \\
\mu_{c2}&=&0.013516 \label{eq:muc2class}\\
\mu_{c3}&=&0.010914\label{eq:muc3class}
\end{eqnarray}
and the $D$ occurring in the second condition of KAM theorem we have
found from (~\ref{eq:D})
\begin{equation}D=\frac{644u^4-541u^2+36}{8(4u^2-1)(25u^2-4)}\label{eq:Dclass}\end{equation}where
$u=\omega_1\omega_2$

All the above results, are exactly similar with the results as in
\cite{Deprit1967}.

\begin{figure}
\epsscale{1}\plotone{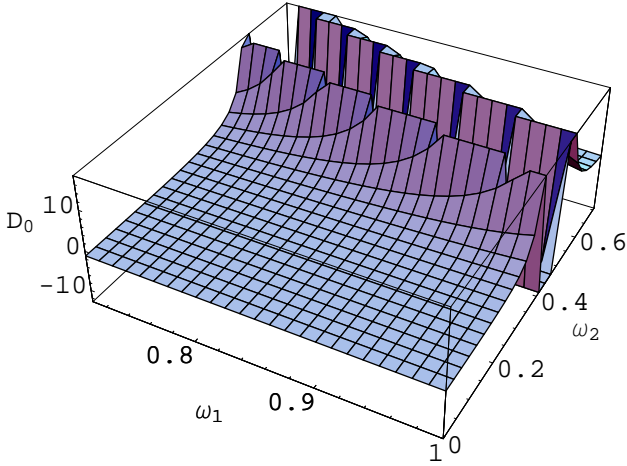} \caption{$A_2=0, q_1=1$}\label{fig:f1}
\end{figure}

\begin{figure}
\plotone{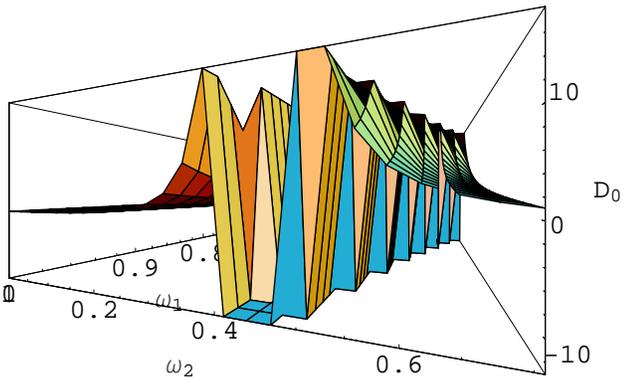} \epsscale{1}\caption{$A_2=0, q_1=1,
\omega_1=0.924270, \omega_2=0.381742,\ D_0= 0$}\label{fig:f2}
\end{figure}

Now we have $\epsilon=1-q_1$,$ W_1=\frac{(1-\mu)\epsilon}{c_d}$,
suppose $D=D_0$. We draw the figure (~\ref{fig:f1})  which describes
the instability range in classical case and figure  (~\ref{fig:f2})
views the points $\omega_1=0.924270,\omega_2=0.381742,D_0= 0$ , when
$A_2=0, q_1=1$ and $\frac{1}{\sqrt{2}}<\omega_1<1,\quad
0<\omega_2<\frac{1}{\sqrt{2}}$ the value of $\gamma=0.978173$.
\subsection{Observation II}\label{subsect:obs2}
\begin{figure}[t]
\epsscale{1} \plotone{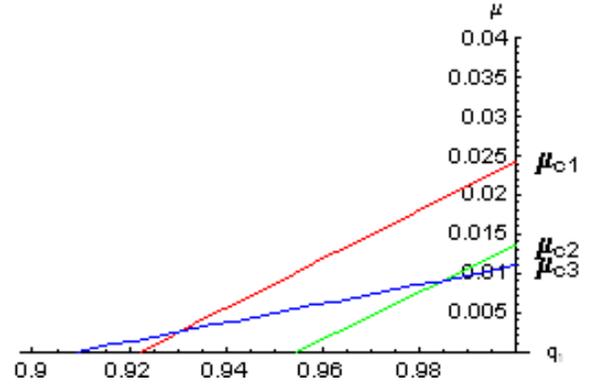} \caption{Stability region
$\mu_{ci},(i=1,2,3)-q_1$ , when $A_2=0$, $W_1\neq0$ }\label{fig:f3}
\end{figure}

 Consider the case when $A_2=0$, $q_1\neq1,(W_1\neq0)$ i.e.
photogravitational restricted three body problem with P-R drag when bigger primary
is supposed to be radiating body and small primary is being
spherical symmetric. The coordinates of triangular equilibrium
points are given by
\begin{equation}
x=x_0\biggl\{1-\displaystyle{\frac{W_1[\amc+\mu\adc]}{3\mu\amc x_0
y_0}}\biggr\}\end{equation}
\begin{equation}
y=y_0\biggl\{1-\displaystyle{\frac{W_1\delta^2[2\mu-1-\mu\adc]}{6\mu\amc
y^3_0}}\biggr\}
\end{equation}
this result coincides with  \citet*{Schuerman1980}, where \(
x=x_0=\adc-\mu,\,
y=y_0=\pm\delta\biggl(1-\frac{\delta^4}{4}\biggr)^{1/2}\),\, $
q_1=1=\delta$,\, \( x=\frac{1}{2}-\mu,\, y=\pm\frac{\sqrt{3}}{2} \)
Substituting $\epsilon=1-q_1$,$ W_1=\frac{(1-\mu)\epsilon}{c_d}$,
$\mu=\mu_{ci},(i=0,1,2,3)$,$A_2=0$  in (~\ref{eq:muc1}),
(~\ref{eq:muc2}) (~\ref{eq:muc3}), we have found that the triangular
equilibrium points are stable  in the range of stability except
three mass ratios
\begin{eqnarray}
\mu_{c1}&=&
 0.024294 -0.312692(1 - q_1)\nonumber\\&& + \frac{0.976732(1 - q_1)}{c_d}\label{eq:muc1pr}\\
\mu_{c2}&=&0.013516 -
    0.29724(1 - q_1)\nonumber\\&& + \frac{0.994062(1 - q_1)}{c_d}\label{eq:muc2pr}\\
\mu_{c3}&=&0.010914-0.120489(1 - q_1)\nonumber\\&&+\frac{2.87259 (1
- q_1)}{c_d}\label{eq:muc3pr}
\end{eqnarray}
We have  observed from table (\ref{tbl-1}) and figure
(~\ref{fig:f3}), the  mass ratio increases,  accordingly as  the  radiation pressure increases,
 these results are similar but not identical
to those of \cite{Papadakis99}.
\subsection{Observation III}\label{subsect:obs3}
\begin{figure}[t]
\epsscale{1} \plotone{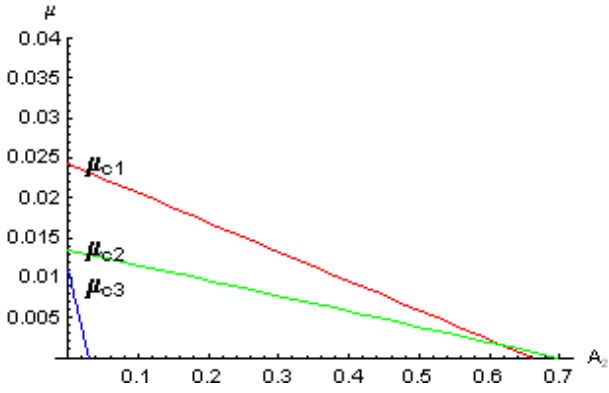} \caption{Stability region
$\mu_{ci},(i=1,2,3)-A_2$, when $ q_1=1$,$W_1=0$}\label{fig:f4}
\end{figure}

When $A_2\neq0$, $q_1=1,(W_1=0)$ i.e. in this observation we have
considered the smaller primary as an oblate spheroid, the radiation
pressure(P-R drag) is not considered. The triangular equilibrium
points are given by
\begin{eqnarray}
x&=&\frac{1-2\mu-A_2}{2}\label{lxob}\\
y&=\pm&\frac{\sqrt{3}}{2}\Bigl\{1-\frac{A_2}{3}\Bigr\}\label{lyob}
\end{eqnarray}
 which are similar but not identical to results as in \cite{Bhatnagaretal1983} and
 \cite{Chandraetal2004}. In this case triangular equilibrium
 points are stable in the nonlinear sense except three mass ratios
 at which Moser's condition fails. Which are given by
 \begin{eqnarray}
\mu_{c1}&=&
 0.024294
   - 0.036851A_2 \\
\mu_{c2}&=&0.013516 - 0.019383 A_2 \\
\mu_{c3}&=&0.010914-0.373118A_2
\end{eqnarray}
The stability region are shown in the diagram
$A_2-\mu_{ci}(i=1,2,3)$, (~\ref{fig:f4}), the outer line is
corresponding to $\mu_{c1}$, second line due to $ \mu_{c2}$ and
innermost line is due to $\mu_{c3}$ it is clear from table
(\ref{tbl-2}) the $\mu$ decreases as $A_2$ increases. These
results agree with \cite{Markellosetal96, Bhatnagaretal1983}
\subsection{Observation IV}\label{subsect:obs4} When $A_2\neq 0,\, q_1\neq1
(W_1\neq0)$ this is the most generalized  case which is being
considered.The triangular equilibrium points are given by
(~\ref{eq:1x}), (~\ref{eq:Ly}) clearly they are the functions of
oblateness coefficient $A_2$ and P-R drag term $W_1$.

\begin{figure*}[h]
\epsscale{1.5} \plottwo{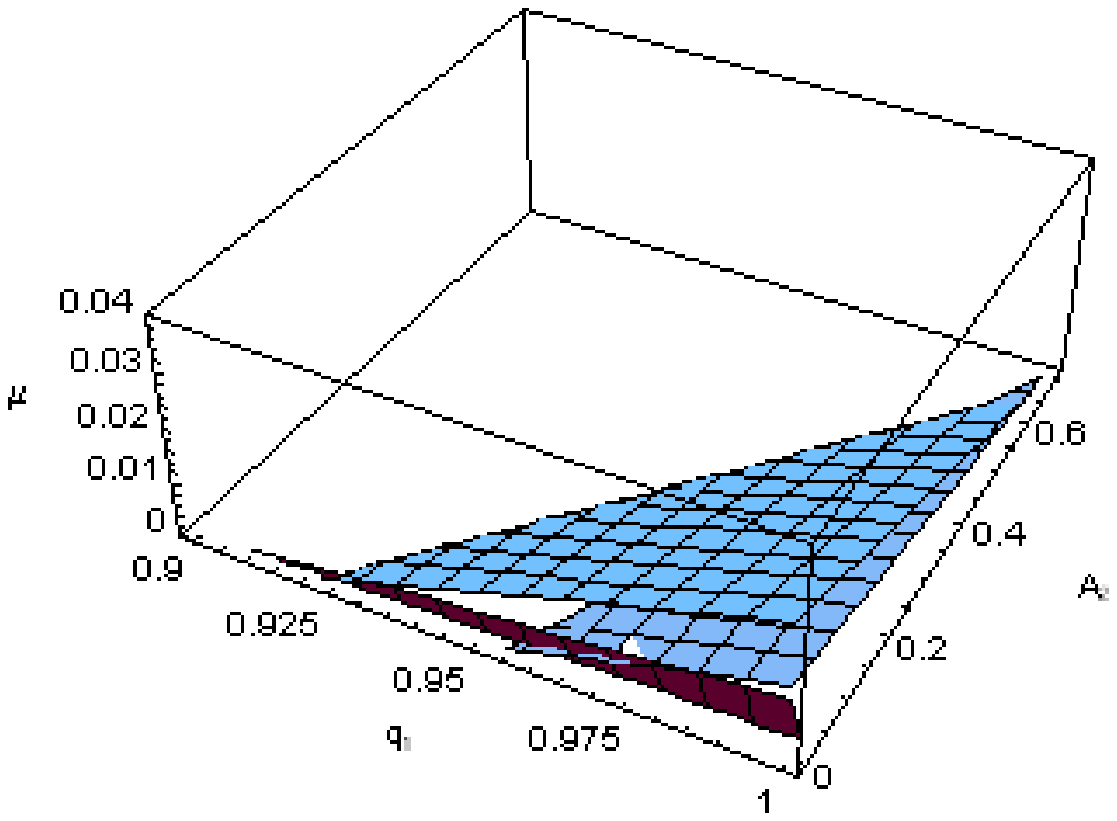}{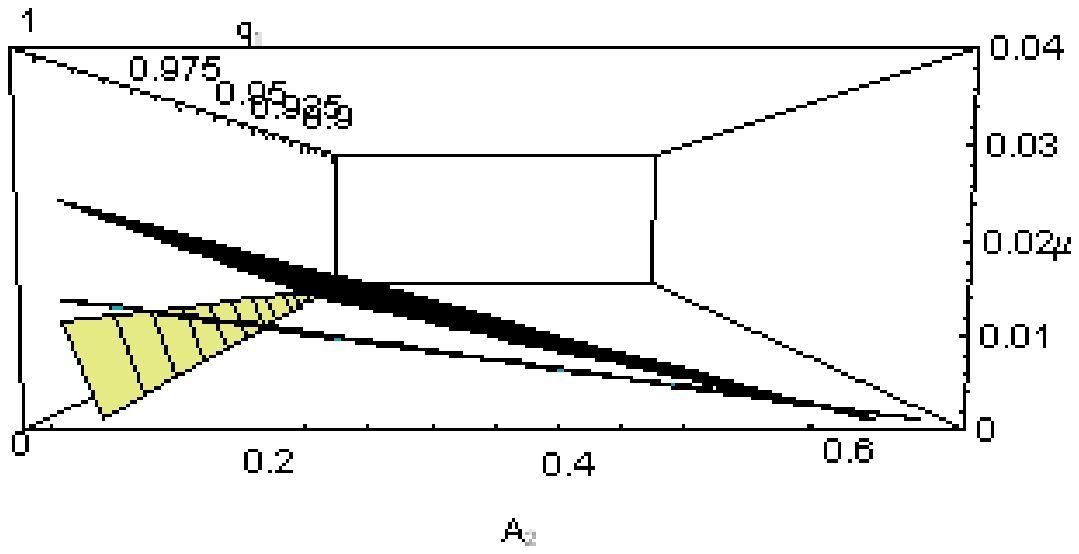} \caption{Both the graphs
 show the stability region $\mu_{ci},(i=1,2,3)-q_1-A_2$, second graph view
 $\mu_{c0}=.035829$}\label{fig:f5}
\end{figure*}

Substituting $\epsilon=1-q_1$,$ W_1=\frac{(1-\mu)\epsilon}{c_d}$,
$\mu=\mu_{ci},(i=0,1,2,3)$ in (~\ref{eq:muc1}), (~\ref{eq:muc2})
(~\ref{eq:muc3}), we get the new formulae
\begin{eqnarray}\mu_{c1}&=&
 0.024294 -
    0.312692(1-q_1)\nonumber\\&-& 0.036851A_2 + \frac{0.976732(1 - q_1)}{c_d}\label{eq:muc1apr}\\
\mu_{c2}&=&0.013516 -
    0.29724(1-q_1)\nonumber\\&-& 0.019383 A_2 + \frac{0.994062(1 - q_1)}{c_d}\label{eq:muc2apr}\\
\mu_{c3}&=&0.010914-0.120489
(1-q_1)\nonumber\\&-&0.373118A_2+\frac{2.87259 (1 -
q_1)}{c_d}\label{eq:muc3apr}\end{eqnarray} Using
(~\ref{eq:muc1apr})-(~\ref{eq:muc3apr}) we have drawn $\mu-A_2-q_1$,
$3D$ diagrams (~\ref{fig:f5}). You can see in the first diagram, the
uppermost plane is due to   $\mu_{c1}$,  middle plane is due to
$\mu_{c2}$ and innermost plane is due to $\mu_{c3}$, second view
value of $\mu_{c0}=.035829$. From these diagrams, we reached at the
conclusion that the stability region is reduced due to P-R drag and
oblateness effect of smaller primary. But still the triangular
equilibrium points are stable in the range of linear stability
except  three mass ratios at which KAM theorem fails, while they are
unstable in linear case [see \cite{Murray1994, KushvahBR2006}].

\section{Conclusion}
\label{sec:conc} Using \cite{Whittaker1965} method we have seen that
the second order part  $H_2$ of the Hamiltonian is transformed into
the normal form $H_2=\omega_1I_1-\omega_2I_2$ and the third order
part $H_3$ of the Hamiltonian in $I_1^{1/2},I_2^{1/2}$  zero. We
conclude that  the stability region is reduced due to P-R drag and
oblateness effect of smaller primary. But still the triangular
equilibrium points are stable in the nonlinear sense in the range of
linear stability except for three mass ratios $\mu_{ci},(i=1,2,3)$
at which KAM theorem fails, while they are unstable in linear case
[see \cite{Murray1994, KushvahBR2006}]. These results  agree with
those found by
 \cite{Deprit1967} and others.

\acknowledgments{ We are thankful to D.S.T. Government of India, New
Delhi for sanctioning a project DST/MS/140/2K dated 02/01/2004 on
this topic. We are also thankful to IUCAA Pune for providing
financial assistance for visiting library and  computer facility.}

\begin{table*}
\tabletypesize{\scriptsize}\caption{$A_2=0,q_1\neq1,(W_1\neq0)$}\label{tbl-1}
\begin{tabular}{crrr}
\tableline\tableline $q_1$ & $\mu_{c1}$ & $\mu_{c2}$  & $\mu_{c3}$ \\
\tableline 0.95 & 0.00866& -0.001346 &
0.00488921 \\
0.96 &0.011786 & 0.0016263 & 0.006094 \\
0.97 & 0.014913 & 0.0045987 &  0.007299\\
0.98 & 0.018040 & 0.0075712&  0.008504\\
0.99 & 0.02117 & 0.010544& 0.00970878 \\
1.00 & 0.024294 & 0.013516 & 0.0109137\\
\tableline
\end{tabular}
\end{table*}
\begin{table*}
\tabletypesize{\scriptsize} \caption{$A_2\neq 0,
q_1=1,(W_1=0)$\label{tbl-2}}
\begin{tabular}{crrr}
\tableline\tableline $A_2$ & $\mu_{c1}$ & $\mu_{c2}$  & $\mu_{c3}$ \\
\tableline
0.0&0.024294 &0.01352 &0.010914\\
0.1 &0.020609&0.01158 &-0.026398\\
0.2&0.016924&0.009639 &-0.06371\\
0.3&0.013239 &0.007701&-0.101022\\
0.4 & 0.009554 &0.005763&-0.138334\\
0.5 &0.005869 &0.003825 &-0.175645\\
0.6&0.002184&0.001886&-0.212957\\
0.7&-0.001501 &-0.000052 &-0.250269\\
\tableline
\end{tabular}
\end{table*}
\appendix
\section*{Appendix I}
Coefficients $A_{1,i},B_{1,i}$ and $C_{1,i},(i=1\dots,7)$ are given
by:
\begin{eqnarray*}
&&A_{1,1}=-\frac{9}{8(-1+2\omega_1^2)^2(-1+5\omega_1^2)}+\frac{259\omega_1^2}{24(-1+2\omega_1^2)^2(-1+5\omega_1^2)}\nonumber\\&&-\frac{205\omega_1^4}{18(-1+2\omega_1^2)^2(-1+5\omega_1^2)}+\frac{31\omega_1^6}{18(-1+2\omega_1^2)^2(-1+5\omega_1^2)}
\end{eqnarray*}
\begin{eqnarray*}
&&A_{1,2}=\frac{1}{36(1-2\omega_1^2)^2(-1+5\omega_1^2)}-\frac{13\omega_1^2}{18(1-2\omega_1^2)^2(-1+5\omega_1^2)}\nonumber\\&&+\frac{13\omega_1^4}{27(1-2\omega_1^2)^2(-1+5\omega_1^2)}+\frac{167\omega_1^6}{72(1-2\omega_1^2)^2(-1+5\omega_1^2)}\nonumber\\&&+\frac{107\omega_1^8}{108(1-2\omega_1^2)^2(-1+5\omega_1^2)}
\end{eqnarray*}
\begin{eqnarray*}
&&A_{1,3}=\frac{1}{2(1-2\omega_1^2)^3(-1+5\omega_1^2)^2}-\frac{421\omega_1^2}{32(1-2\omega_1^2)^3(-1+5\omega_1^2)^2}\nonumber\\&&-\frac{19\omega_1^4}{2(1-2\omega_1^2)^3(-1+5\omega_1^2)^2}-\frac{8141559\omega_1^6}{32(1-2\omega_1^2)^3(-1+5\omega_1^2)^2}\nonumber\\&&+\frac{29\omega_1^8}{(1-2\omega_1^2)^3(-1+5\omega_1^2)^2}
\end{eqnarray*}
\begin{eqnarray*}
A_{1,4}&=&\frac{1319}{436(1-2\omega_1^2)^2(-1+5\omega_1^2)}-\frac{12639\omega_1^2}{436(1-2\omega_1^2)^2(-1+5\omega_1^2)}\nonumber\\&&+\frac{14275\omega_1^4}{436(1-2\omega_1^2)^2(-1+5\omega_1^2)}-\frac{799\omega_1^6}{218(1-2\omega_1^2)^2(-1+5\omega_1^2)}
\end{eqnarray*}
\begin{eqnarray*}
A_{1,5}&=&\frac{57}{52(1-2\omega_1^2)^3(-1+5\omega_1^2)^2}-\frac{525\omega_1^2}{52(1-2\omega_1^2)^3(-1+5\omega_1^2)^2}\nonumber\\&&-\frac{475\omega_1^4}{26(1-2\omega_1^2)^3(-1+5\omega_1^2)^2}+\frac{1559\omega_1^6}{26(1-2\omega_1^2)^3(-1+5\omega_1^2)^2}\nonumber\\&&+\frac{283\omega_1^8}{13(1-2\omega_1^2)^3(-1+5\omega_1^2)^2}
\end{eqnarray*}
\begin{eqnarray*}
A_{1,6}&=&-\frac{2747\omega_1^2}{10368\sqrt{3}(-1+2\omega_1^2)}+\frac{41(9+4\omega_1^2)}{9216\sqrt{3}(-1+2\omega_1^2)^2}\nonumber\\&&-\frac{93899(9+4\omega_1^2)}{331776\sqrt{3}(-1+2\omega_1^2)}+\frac{12875\omega_1^2(9+4\omega_1^2)}{82944\sqrt{3}(-1+2\omega_1^2)^2}
\end{eqnarray*}
\begin{eqnarray*}
A_{1,7}&=&-\frac{1337}{6144\sqrt{3}(-1+2\omega_1^2)}+\frac{779\omega_1(9+4\omega_1^2)}{10368\sqrt{3}(-1+2\omega_1^2)^2}+\frac{41(9+4\omega_1^2)}{18432\sqrt{3}(-1+2\omega_1^2)^2}\nonumber\\&&-\frac{227347\omega_1^2(9+4\omega_1^2)}{331776\sqrt{3}(-1+2\omega_1^2)}-\frac{37259(9+4\omega_1^2)}{82944\sqrt{3}(-1+2\omega_1^2)}\nonumber\\&&+\frac{6517\omega_1^2(9+4\omega_1^2)}{3072\sqrt{3}(-1+2\omega_1^2)^2(4\omega_1^2-\omega_2^2)}
\end{eqnarray*}
\begin{eqnarray*}
B_{1,1}&=&\frac{43\omega_1\omega_2}{6(1-5\omega_1^2)(-1+2\omega_1^2)(1-5\omega_2^2)(1-2\omega_2^2)}\nonumber\\&&+\frac{32\omega_1^3\omega_2^3}{3((1-5\omega_1^2))(-1+2\omega_1^2)(1-5\omega_2^2)(1-2\omega_2^2)}
\end{eqnarray*}
\begin{eqnarray*}
B_{1,2}&=&\frac{309\omega_1\omega_2}{8(1-5\omega_1^2)(1-2\omega_1^2)(1-5\omega_2^2)(1-2\omega_2^2)}\nonumber\\&&+\frac{5904\omega_1\omega_2}{(-1+2\omega_1^2)(9+4\omega_2^2)^2}-\frac{407\omega_1^3\omega_2^3}{6(1-5\omega_1^2)(1-2\omega_1^2)(1-5\omega_2^2)(1-2\omega_2^2)}
\end{eqnarray*}
\begin{eqnarray*}
B_{1,3}&=&\frac{1800\omega_1\omega_2}{(-1+2\omega_1^2)(9+4\omega_2^2)^2}\nonumber\\&&+\frac{10083-614070\omega_1^2\omega_2^2+400800\omega_1^4\omega_2^4-3035216\omega_1^6\omega_2^6-260802\omega_1^8\omega_2^8}{8\omega_1\omega_2(9-59\omega_1^2+62\omega_1^4+40\omega_1^6)(9-59\omega_2^2+62\omega_2^4+40\omega_2^6)}
\end{eqnarray*}
\begin{eqnarray*}
B_{1,4}&=&\frac{247\omega_1\omega_2}{4(1-5\omega_1^2)(1-2\omega_1^2)(1-5\omega_2^2)(1-2\omega_2^2)}\nonumber\\&&+\frac{6817\omega_1^3\omega_2^3}{36(1-5\omega_1^2)(1-2\omega_1^2)(1-5\omega_2^2)(1-2\omega_2^2)}
\end{eqnarray*}
\begin{eqnarray*}
B_{1,5}&=&\frac{1800\omega_1\omega_2}{(-1+2\omega_1^2)(9+4\omega_2^2)^2}\nonumber\\&&+\frac{-89211+2042998\omega_1^2\omega_2^2+1028577\omega_1^4\omega_2^4\omega_1^2+16052098\omega_1^6\omega_2^6+1215804\omega_1^8\omega_2^8}{32\omega_1\omega_2(-1+5\omega_1^2)^2(-1+5\omega_2^2)^2(-9+14\omega_1^2+8\omega_1^4)(-9-14\omega_2^2+8\omega_2^4)}
\end{eqnarray*}
\begin{eqnarray*}
B_{1,6}&=&\frac{1599\sqrt{3}(9+192\omega_1\omega_2+\omega_2^2)}{512\omega_1\omega_2(-1+2\omega_1^2)(9+4\omega_2^2)^2}\nonumber
\end{eqnarray*}
\begin{eqnarray*}
B_{1,7}&=-&\frac{3\sqrt{3}(2398599-9031680\omega_2^2-369\omega_1\omega_2^3+574\omega_1\omega_2^5+15744\omega_1^2\omega_2^6+328\omega_2^7}{512\omega_1^2\omega_2^4(-1+2\omega_1^2)(-9+14\omega_2^2+8\omega_2^4)}\nonumber\\&&-\frac{192(-41601+41\omega_1^2)\omega_2^4}{512\omega_1^2\omega_2^4(-1+2\omega_1^2)(-9+14\omega_2^2+8\omega_2^4)}
\end{eqnarray*}
\begin{eqnarray*}
C_{1,1}&=&\frac{9}{8(-1+2\omega_2^2)^2(-1+5\omega_2^2)}+\frac{205\omega_2^2}{24(-1+2\omega_2^2)^2(-1+5\omega_2^2)}\nonumber\\&&-\frac{205\omega_2^4}{18(-1+2\omega_2^2)^2(-1+5\omega_2^2)}+\frac{31\omega_2^6}{18(-1+2\omega_1^2)^2(-1+5\omega_1^2)}
\end{eqnarray*}
\begin{eqnarray*}
C_{1,2}&=&\frac{1}{36(1-2\omega_2^2)^2(-1+5\omega_2^2)}-\frac{13\omega_2^2}{18(1-2\omega_2^2)^2(-1+5\omega_2^2)}\nonumber\\&&+\frac{13\omega_2^4}{27(1-2\omega_2^2)^2(-1+5\omega_2^2)}-\frac{167\omega_2^6}{72(1-2\omega_2^2)^2(-1+5\omega_2^2)}\nonumber\\&&+\frac{107\omega_2^8}{108(1-2\omega_2^2)^2(-1+5\omega_2^2)}
\end{eqnarray*}
\begin{eqnarray*}
C_{1,3}&=&\frac{1}{2(1-2\omega_2^2)^3(-1+5\omega_2^2)^2}-\frac{421\omega_2^2}{32(1-2\omega_2^2)^3(-1+5\omega_2^2)^2}\nonumber\\&&-\frac{19\omega_2^4}{2(1-2\omega_2^2)^3(-1+5\omega_2^2)^2}-\frac{407\omega_2^6}{16(1-2\omega_2^2)^3(-1+5\omega_2^2)^2}\nonumber\\&&+\frac{29\omega_2^8}{(1-2\omega_1^2)^3(-1+5\omega_1^2)^2}
\end{eqnarray*}
\begin{eqnarray*}
C_{1,4}&=&\frac{1319}{436(1-2\omega_2^2)^2(-1+5\omega_2^2)}-\frac{12639\omega_2^2}{436(1-2\omega_2^2)^2(-1+5\omega_2^2)}\nonumber\\&&+\frac{14275\omega_2^4}{436(1-2\omega_2^2)^2(-1+5\omega_2^2)}-\frac{799\omega_2^6}{218(1-2\omega_2^2)^2(-1+5\omega_2^2)}
\end{eqnarray*}
\begin{eqnarray*}
C_{1,5}&=&\frac{57}{52(1-2\omega_2^2)^3(-1+5\omega_2^2)^2}+\frac{525\omega_2^2}{52(1-2\omega_2^2)^3(-1+5\omega_2^2)^2}\nonumber\\&&-\frac{475\omega_2^4}{26(1-2\omega_2^2)^3(-1+5\omega_2^2)^2}+\frac{1559\omega_2^6}{26(1-2\omega_2^2)^3(-1+5\omega_2^2)^2}\nonumber\\&&+\frac{283\omega_2^8}{13(1-2\omega_2^2)^3(-1+5\omega_2^2)^2}
\end{eqnarray*}
\begin{eqnarray*}
C_{1,6}&=&-\frac{287\sqrt{3}(-3+32\omega_2^2+48\omega_2^4)}{1024\omega_2^2(-9+14\omega_2^2+8\omega_2^4)}\nonumber
\end{eqnarray*}
\begin{eqnarray*}
C_{1,7}&=-&\frac{\sqrt{3}82\omega_1^2(3-38\omega_2^2+16\omega_2^4+96\omega_2^6)}{512(9+4\omega_2^2)(-\omega_1^2+4\omega_2^2)(\omega_2^2-2\omega_2^3)^2}\nonumber\\&&+\frac{3\sqrt{3}\omega_2^2(-142911+195110\omega_2^2+74728\omega_2^4+66784\omega_2^6)}{512(9+4\omega_2^2)(-\omega_1^2+4\omega_2^2)(\omega_2^2-2\omega_2^3)^2}
\end{eqnarray*}
\section*{Appendix II}
\begin{eqnarray*}
D_2&=&\frac{567(-151+16\omega_1^2)}{16384(-1+2\omega_1^2)^2(9+4\omega_1^2)^2}+\frac{\omega_1^2\omega_2^2}{884736}\biggl\{\frac{1620864}{(-1+2\omega_1^2)^2}+\frac{2507364}{(1-2\omega_1^2)}\\&&+\frac{706482}{(-1+2\omega_1^2)^2(4\omega_1^2-\omega_2^2)}+\frac{71663616000}{(-1+2\omega_1^2)(9+4\omega_2^2)}\nonumber\\&&+\frac{8062156800}{(-1+2\omega_2^2)(9+4\omega_2^2)^2}+\frac{1074954240}{(-1+2\omega_1^2)^2(9+4\omega_2^2)}\\&&+\frac{112969617408}{(1-5\omega_1^2)^2(-1+2\omega_1^2)^2(9+4\omega_1^2)(1-5\omega_2^2)^2(-1+2\omega_2^2)^2(9+4\omega_2^2)}\\&&+\frac{17146183680}{(9-14\omega_2^2-8\omega_2^4)}\biggr\}\\&&+\frac{1028577\omega_1^4\omega_2^4}{16(1-5\omega_1^2)^2(-1+2\omega_1^2)^2(9+4\omega_1^2)(1-5\omega_2^2)^2(-1+2\omega_2^2)^2(9+4\omega_2^2)}\nonumber\\&&+\frac{8026049\omega_1^6\omega_2^6}{8(1-5\omega_1^2)^2(-1+2\omega_1^2)^2(9+4\omega_1^2)(1-5\omega_2^2)^2(-1+2\omega_2^2)^2(9+4\omega_2^2)}\nonumber\\&&+\frac{303951\omega_1^8\omega_2^8}{4(1-5\omega_1^2)^2(-1+2\omega_1^2)^2(9+4\omega_1^2)(1-5\omega_2^2)^2(-1+2\omega_2^2)^2(9+4\omega_2^2)}
\end{eqnarray*}
\begin{eqnarray*}
D_3&=&\frac{3}{8192(-1+2\omega_1^2)}\Bigl\{819+\frac{8064}{(2\omega_1+\omega_2)(\omega_1+2\omega_2)(9+4\omega_2)^2}\nonumber\\&&-\frac{6883328}{(1-5\omega_1^2)(-1+2\omega_1^2)(9+4\omega_1^2)(1-5\omega_2^2)(-1+2\omega_2^2)(9+4\omega_2^2)}
\Bigr\}\nonumber\\&&+\frac{\omega_1^2\omega_2^2}{147456}\Bigl\{\frac{706240}{(-1+2\omega_1^2)}+\frac{289737}{(1-2\omega_1^2)}(4\omega_1^2-\omega_2^2)-\frac{530841600}{(-1+2\omega_1^2)(9+4\omega_2^2)^2}\nonumber\\&&+\frac{59719680}{(-1+2\omega_2^2)(9+4\omega_2^2)^2}+\frac{59719680\omega_2^2}{(-1+2\omega_1^2)(9+4\omega_2^2)^2}+\frac{3317760}{(-1+2\omega_2^2)^2(9+4\omega_2^2)}\\&&+\frac{71516160}{9-14\omega_2^2-8\omega_2^4}+\frac{24772608}{(\omega_1^2-4\omega_2^2)(-1+2\omega_2^2)(9+4\omega_2^2)}\\&&+\frac{22637076480}{(1-5\omega_1^2)(-1+2\omega_1^2)(9+4\omega_1^2)(1-5\omega_2^2)(-1+2\omega_2^2)(9+4\omega_2^2)}\Bigr\}\nonumber\\&&-\frac{100200\omega_1^4\omega_2^4}{(1-5\omega_1^2)(-1+2\omega_1^2)(9+4\omega_1^2)(1-5\omega_2^2)(-1+2\omega_2^2)(9+4\omega_2^2)}\nonumber\\&&+\frac{758804\omega_1^6\omega_2^6}{(1-5\omega_1^2)(-1+2\omega_1^2)(9+4\omega_1^2)(1-5\omega_2^2)(-1+2\omega_2^2)(9+4\omega_2^2)}\nonumber\\&&+\frac{130401\omega_1^8\omega_2^8}{2(1-5\omega_1^2)(-1+2\omega_1^2)(9+4\omega_1^2)(1-5\omega_2^2)(-1+2\omega_2^2)(9+4\omega_2^2)}\end{eqnarray*}\begin{eqnarray*}
D_4&=&\frac{1}{294912}\Biggl\{243\biggl\{\frac{58477}{(1-2\omega_1^2)^2}+\frac{89216}{(9-14\omega_1^2-8\omega_1^4)}+\frac{7872}{(-1+2\omega_2^2)^2}\\&&+\frac{33456}{(9-14\omega_2^2-8\omega_2^4)}\biggr\}\\&&+2\omega_1^2\omega_2^2\biggl\{\frac{5864788}{(1-2\omega_1^2)}-\frac{186165}{(-1+2\omega_1^2)^2(4\omega_1^2-\omega_2^2)}+\frac{1885814784}{(\omega_1^2-4\omega_2^2)(9-14\omega_2^2-8\omega_2^4)}\\&&+\frac{18210816}{(1-7\omega_1^2+10\omega_1^4)(1-7\omega_2^2+10\omega_2^4)}\biggr\}-\frac{111689728\omega_1^4\omega_2^4}{(1-7\omega_1^2+10\omega_1^4)(1-7\omega_2^2+10\omega_2^4)}\Bigg\}
\end{eqnarray*}
\begin{eqnarray*}
D_5&=&\frac{1}{49152}\Biggl\{9\biggl\{-\frac{2457}{(1-2\omega_1^2)^2}+\frac{6426}{(-9+14\omega_1^2+8\omega_1^4)}\\&&-\frac{30450688}{(-1+5\omega_1^2)^2(9-14\omega_1^2-8\omega_1^4)(-1+5\omega_2^2)^2(9-14\omega_2^2-8\omega_2^4)}\biggr\}\\&&+\omega_1^2\omega_2^2\biggl\{\frac{90048}{(1-2\omega_1^2)^2}+\frac{139298}{(1-2\omega_1^2)^2}+\frac{39249}{(1-2\omega_1^2)^2(4\omega_1^2-\omega_2^2)}\\&&+\frac{447897600}{(-1+2\omega_1^2(9+4\omega_2^2)}+\frac{952565760}{(9-14\omega_2^2-8\omega_2^4)}\\&&+\frac{6276089856}{(1-5\omega_1^2)^2(-9+14\omega_1^2+8\omega_1^4)(1-5\omega_2^2)^2(-9+14\omega_2^2+8\omega_2^4)}\\&&-\frac{594542592}{(\omega_1^2-4\omega_2^2)(-9+14\omega_2^2+8\omega_2^4)}\biggr\}\\&&+\frac{3159788544\omega_1^4\omega_2^4}{(1-5\omega_1^2)^2(-9+14\omega_1^2+8\omega_1^4)(1-5\omega_2^2)^2(-9+14\omega_2^2+8\omega_2^4)}\\&&+\frac{49312045056\omega_1^6\omega_2^6}{(1-5\omega_1^2)^2(-9+14\omega_1^2+8\omega_1^4)(1-5\omega_2^2)^2(-9+14\omega_2^2+8\omega_2^4)}\\&&+\frac{3734949888\omega_1^8\omega_2^8}{(1-5\omega_1^2)^2(-9+14\omega_1^2+8\omega_1^4)(1-5\omega_2^2)^2(-9+14\omega_2^2+8\omega_2^4)}\Bigg\}
\end{eqnarray*}
\begin{eqnarray*}
D_6&=&\frac{1}{82944\sqrt{3}}\Biggl\{29889\biggl\{\frac{52}{(-1+2\omega_1^2)^2}+\frac{7}{(9-14\omega_2^2-8\omega_2^4)}\biggr\}\\&&++2\omega_1^2\omega_2^2\biggl\{-\frac{738}{(1-2\omega_1^2)^2}+\frac{93899}{(-1+2\omega_1^2)^2}+\frac{91445760}{(\omega_1^2-4\omega_2^2)(9-14\omega_2^2-8\omega_2^4)}\Bigg\}
\end{eqnarray*}
\begin{eqnarray*}
D_7&=&\frac{1}{110592\sqrt{3}}\Biggl\{27\biggl\{\frac{5904}{(-1+2\omega_1^2)^2}+\frac{122157}{(\omega_1^2-4\omega_2^2)(-1+2\omega_1^2)^2}\\&&-\frac{758086}{(\omega_1^2-4\omega_2^2)(-1+2\omega_1^2)}+\frac{5904}{(9-14\omega_2^2-8\omega_2^4)}\biggr\}\\&&+2\omega_1^2\omega_2^2\biggl\{-\frac{492}{(1-2\omega_1^2)^2}+\frac{370964}{(-1+2\omega_1^2)}-\frac{58653}{(4\omega_1^2-\omega_2^2)(-1+2\omega_1^2)^2}\\&&+\frac{13893120}{(9+4\omega_2^2)^2(1-2\omega_2^2)}-\frac{116702208}{(\omega_1^2-4\omega_2^2)(-1+2\omega_2^2)^2(9+4\omega_2^2)^2}\\&&+\frac{103680}{(9+4\omega_2^2)(1-2\omega_2^2)^2}+\frac{870912}{(\omega_1^2-4\omega_2^2)(-1+2\omega_2^2)^2(9+4\omega_2^2)^2}\\&&+\frac{62519040}{(-9+14\omega_2^2+8\omega_2^4)}-\frac{246177792}{(\omega_1^2-4\omega_2^2)(9+4\omega_2^2)^2}\biggr\}\Bigg\}
\end{eqnarray*}

\end{document}